\documentclass[MSNbibl,number,citesort,dvips]{arxbj}
\usepackage{upgreek}

\aid{0}
\volume{19}
\issue{5A}
\pubyear{2013}
\firstpage{2000}
\lastpage{2009}
\doi{10.3150/12-BEJ438} 

\makeatletter
\newcommand{\rrvert}{\vert}
\newcommand{\llvert}{\vert}
\newtheorem{prop}{Proposition}[section]
\newremark{rem}[prop]{Remark}
\newtheorem{theo}[prop]{Theorem}
\newtheorem{cor}[prop]{Corollary}
\makeatother

\begin{document}
\begin{frontmatter}

\title{Integrability properties and limit theorems for the exit time
from a cone of planar Brownian motion}
\runtitle{Integrability properties and limit theorems}

\begin{aug}
\author[1,2]{\fnms{Stavros} \snm{Vakeroudis}\corref{}\thanksref{1,2,e1}\corref{}\ead[label=e1,mark]{stavros.vakeroudis@upmc.fr}\ead[label=u1,url]{http://svakeroudis.wordpress.com}}
\and
\author[1,3]{\fnms{Marc} \snm{Yor}\thanksref{1,3,e2}\ead[label=e2,mark]{yormarc@aol.com}}
\runauthor{S. Vakeroudis and M. Yor} 
\address[1]{Laboratoire de Probabilit\'{e}s et Mod\`{e}les
Al\'{e}atoires (LPMA), CNRS: UMR7599, Universit\'{e} Pierre et Marie
Curie -- Paris VI, Universit\'{e} Paris-Diderot Paris VII,
4, Place Jussieu, 75252 Paris Cedex 05, France. \mbox{\printead
{e1,u1}}; \printead*{e2}}
\address[2]{Probability and Statistics Group, School of Mathematics,
University of Manchester,
Alan Turing Building, Oxford Road, Manchester M13 9PL, United Kingdom}
\address[3]{Institut Universitaire de France, Paris, France}
\end{aug}

\received{\smonth{1} \syear{2012}}

%
\begin{abstract}
We obtain some integrability properties and some limit theorems for
the exit time from a cone of a planar Brownian motion, and we check that
our computations are correct via Bougerol's identity.
\end{abstract}

%

\end{frontmatter}
%

\section{Introduction}\label{sec1}

We consider a standard planar Brownian motion
$(Z_{t}=X_{t}+\mathrm{i}Y_{t},t\geq0)$, starting from \mbox{$x_{0}+\mathrm{i}0,x_{0}>0$}, where
$(X_{t},t\geq0)$ and $(Y_{t},t\geq0)$ are two independent linear
Brownian motions, starting respectively, from $x_{0}$ and $0$
(when we simply write: Brownian motion, we always mean real-valued
Brownian motion, starting from 0;
for 2-dimensional Brownian motion, we indicate planar or complex BM).

As is well known It{\^o} and McKean \cite{ItMK65}, since
$x_{0}\neq0$, $(Z_{t},t\geq0)$ does not visit a.s. the point $0$
but keeps winding around $0$ infinitely often. In particular, the
continuous winding process
$\theta_{t}=\operatorname{Im}(\int^{t}_{0}\frac{\mathrm{d}Z_{s}}{Z_{s}}),t\geq0$
is well defined. A scaling argument shows that we may assume $x_{0}=1$,
without loss of generality, since, with obvious notation:
%
\begin{eqnarray}
\bigl(Z^{(x_{0})}_{t},t\geq0 \bigr)\stackrel{(\mathit{law})} {=}
\bigl(x_{0}Z^{(1)}_{(t/x^{2}_{0})},t\geq0 \bigr).
\end{eqnarray}
Thus, from now on, we shall take $x_{0}=1$.

Furthermore, there is the skew product representation:
%
\begin{eqnarray}
\label{skew-product}
\log\llvert Z_{t}\rrvert +\mathrm{i}\theta_{t}
\equiv\int^{t}_{0}\frac
{\mathrm{d}Z_{s}}{Z_{s}}= (
\beta_{u}+\mathrm{i}\gamma_{u} ) \bigg|_{u=H_{t}=\int^{t}_{0}\frac{\mathrm{d}s}{\llvert Z_{s}\rrvert ^{2}}}  ,
\end{eqnarray}
where $(\beta_{u}+\mathrm{i}\gamma_{u},u\geq0)$ is another planar
Brownian motion starting from $\log1+\mathrm{i}0=0$. Thus, the Bessel clock $H$
plays a key role
in many aspects of the study of the winding number process $(\theta_{t},t\geq0)$ (see, e.g., Yor \cite{Yor80}).

Rewriting (\ref{skew-product}) as:
%
\begin{equation}
\label{skew-product2} \log\llvert Z_{t}\rrvert =\beta_{H_{t}}; \qquad   \theta_{t}=\gamma_{H_{t}} ,
\end{equation}
we easily obtain that the two $\sigma$-fields $\sigma\{\llvert Z_{t}\rrvert ,t\geq0\}$ and $\sigma\{\beta_{u},u\geq0\}$ are
identical, whereas $(\gamma_{u},u\geq0)$ is independent from $(\llvert Z_{t}\rrvert ,t\geq0)$.

We shall also use Bougerol's celebrated identity in law, see, for
example, Bougerol \cite{Bou83}, Alili, Dufresne and Yor \cite{ADY97} and Yor \cite{Yor01} (page 200), which may be
written as:
%
\begin{equation}
\label{boug} \mbox{for~fixed }   t \qquad   \sinh(\beta_{t})
\stackrel {(\mathit{law})} {=} \hat{\beta}_{A_{t}(\beta)},
\end{equation}
where $(\beta_{u},u\geq0)$ is 1-dimensional BM, $A_{u}(\beta)=\int^{u}_{0}\mathrm{d}s \exp(2\beta_{s})$ and $(\hat{\beta}_{v},v\geq0)$ is
another BM, independent of $(\beta_{u},u\geq0)$. For the random times
$T^{|\theta|}_{c}\equiv\inf\{ t\dvt|\theta_{t}|=c \}$, and
$T^{|\gamma|}_{c}\equiv\inf\{ t\dvt|\gamma_{t}|=c \}$, $(c>0)$ by
using the skew-product
representation (\ref{skew-product2}) of planar Brownian motion Revuz and Yor \cite
{ReY99}, we obtain:
%
\begin{eqnarray}
\label{skew-productplanar} T^{|\theta|}_{c}=A_{T^{|\gamma|}_{c}}(\beta)
\equiv\int^{T^{|\gamma
|}_{c}}_{0}\mathrm{d}s\exp(2\beta_{s})=H^{-1}_{u}
\bigg|_{u=T^{|\gamma|}_{c}}  .
\end{eqnarray}
Moreover, it has been recently shown that, Bougerol's identity applied
with the random time $T^{|\theta|}_{c}$
instead of $t$ in (\ref{boug}) yields the following Vakeroudis \cite{Vak11}.

\begin{prop}\label{proppseudo}
The distribution of $T^{|\theta|}_{c}$ is characterized by its
Gauss--Laplace transform:
%
\begin{eqnarray}
\label{pseudoTL} E \biggl[ \sqrt{\frac{2c^{2}}{\uppi T^{|\theta|}_{c}}} \exp \biggl( -
\frac{x}{2T^{|\theta|}_{c}} \biggr) \biggr] = \frac{1}{\sqrt{1+x}} \varphi_{m}(x)
\end{eqnarray}
for every $x \geq0$, with $m=\frac{\uppi}{2c}$, and:
%
\begin{eqnarray}
\varphi_{m}(x)=\frac{2}{(G_{+}(x))^{m}+(G_{-}(x))^{m}}, \qquad    G_{\pm
}(x)=
\sqrt{1+x}\pm\sqrt{x}.
\end{eqnarray}
\end{prop}

The remainder of this article is organized as follows: in Section \ref
{chIPLT}, we study
some integrability properties for the exit times from a cone;
more precisely, we obtain some new results concerning the negative
moments of
$T^{|\theta|}_{c}$ and of $T^{\theta}_{c}\equiv\inf\{ t\dvt\theta_{t}=c \}$.
In Section~\ref{secLT}, we state and prove some Limit theorems for
these random times for
$c\rightarrow0$ and for $c\rightarrow\infty$ followed by several
generalizations
(for extensions of these works to more general planar processes, see,
e.g., Doney and Vakeroudis \cite{DoV12}).
We use these results in order to obtain (see Remark \ref{rem}) a new
simple non-computational
proof of Spitzer's celebrated asymptotic theorem Spitzer \cite{Spi58}, which
states that:
%
\begin{eqnarray}
\label{Spi} \frac{2}{\log t}   \theta_{t}
\mathop{\stackrel{(\mathit{law})}{\longrightarrow}}_{t\rightarrow\infty}
 C_{1}  ,
\end{eqnarray}
with $C_{1}$ denoting a standard Cauchy variable
(for other proofs, see, e.g., Williams \cite{Wil74}, Durrett \cite{Dur82}, Messulam and Yor \cite{MeY82}, Bertoin and Werner \cite{BeW94},
Yor \cite{Yor97}, Vakeroudis \cite{Vak11}).
Finally, in Section \ref{secchecks}, we use the Gauss--Laplace
transform (\ref{pseudoTL}) which is equivalent to
Bougerol's identity (\ref{boug}) in order to check our results.

\section{Integrability properties}\label{chIPLT}

Concerning the moments of $T^{|\theta|}_{c}$, we have the following
(a more extended discussion is found in, e.g., Matsumoto and Yor \cite{MaY05}).

\begin{theo} \label{theomoments}
For every $c>0$, $T^{|\theta|}_{c}$ enjoys the following integrability
properties:
\begin{enumerate}[(ii)]
\item[(i)] for $p>0$, $E [ (T^{|\theta|}_{c} )^{p}
]<\infty$, if and only if $p<\frac{\uppi}{4c}$,
\item[(ii)] for any $p<0$, $E [ (T^{|\theta|}_{c}
)^{p} ]<\infty$.
\end{enumerate}
\end{theo}

\begin{cor} \label{momnegtcd}
For $0<c<d$, the random times $T^{\theta}_{-d,c}\equiv\inf\{
t\dvt\theta_{t}\notin(-d,c) \}$, $T^{|\theta|}_{c}$ and $T^{\theta}_{c}$
satisfy the inequality:
%
\begin{eqnarray}
T^{\theta}_{c}\geq T^{\theta}_{-d,c} \geq
T^{|\theta|}_{c}.
\end{eqnarray}
Thus, their negative moments satisfy:
%
\begin{eqnarray}
\mbox{for }  p>0\qquad   E \biggl[\frac{1}{ (T^{\theta}_{c}
)^{p}} \biggr]\leq E \biggl[\frac{1}{ (T^{\theta}_{-d,c}
)^{p}}
\biggr]\leq E \biggl[\frac{1}{ (T^{|\theta|}_{c}
)^{p}} \biggr]<\infty.
\end{eqnarray}
\end{cor}

\begin{pf*} {Proofs of Theorem \ref{theomoments} and
of Corollary \ref{momnegtcd}}
\begin{longlist}[(ii)]
\item[(i)] The original proof is given by Spitzer \cite{Spi58},
followed later by many authors Williams \cite{Wil74}, Burkholder \cite{Bur77}, Messulam and Yor \cite{MeY82},
Durrett \cite{Dur82}, Yor \cite{Yor85}.
See also Revuz and Yor \cite{ReY99}, Ex.~2.21, page 196.

\item[(ii)] In order to obtain this result, we might use the
representation $T^{|\theta|}_{c}=A_{T^{|\gamma|}_{c}}$ together with
a recurrence
formula for the negative moments of $A_{t}$ \cite{Duf00}, Theorem 4.2,
page 417 (in fact, Dufresne also
considers $A^{(\mu)}_{t}=\int^{t}_{0}\mathrm{d}s \exp(2\beta_{s}+2\mu s)$,
but we only need to take $\mu=0$ for our purpose,
and we note $A_{t}\equiv A^{(0)}_{t}$), \cite{Vakth11}. However, we can
also obtain this result by simply remarking
that the RHS of the Gauss--Laplace transform (\ref{pseudoTL}) in
Proposition \ref{proppseudo} is
an infinitely differentiable function in 0 (see also \cite{VaY11}), thus:
%
\begin{eqnarray}
E \biggl[\frac{1}{ (T^{|\theta|}_{c} )^{p}} \biggr]<\infty \qquad   \mbox{for every }  p>0.
\end{eqnarray}
\end{longlist}
Now, Corollary \ref{momnegtcd} follows immediately from Theorem \ref
{theomoments}(ii).
\end{pf*}

\section{\texorpdfstring{Limit theorems for $T^{|\theta|}_{c}$}
{Limit theorems for T|theta|c}}\label{secLT}

\subsection{\texorpdfstring{Limit theorems for $T^{|\theta|}_{c}$, as $c\rightarrow 0$ and $c\rightarrow\infty$}
{Limit theorems for T|theta|c, as c -> 0$ and c -> infinity}}

The skew-product representation of planar Brownian motion
allows to prove the three following asymptotic results for $T^{|\theta|}_{c}$.

\begin{prop}\label{prop1}
\begin{longlist}[(b)]
\item[(a)] For $c\rightarrow0$, we have:
%
\begin{eqnarray}
\frac{1}{c^{2}}   T^{|\theta|}_{c}
\mathop{\stackrel{(\mathit{law})}{\longrightarrow}}_{c\rightarrow0}
 T^{|\gamma|}_{1}.
\end{eqnarray}
\item[(b)] For $c\rightarrow\infty$, we have:
%
\begin{eqnarray}
\frac{1}{c}   \log \bigl(T^{|\theta|}_{c} \bigr) \mathop{\stackrel{(\mathit{law})}{\longrightarrow}}_{c\rightarrow\infty} 2|\beta
|_{T^{|\gamma|}_{1}}.
\end{eqnarray}
\item[(c)] For $\varepsilon\rightarrow0$, we have:
%
\begin{eqnarray}
\frac{1}{\varepsilon^{2}}   \bigl(T^{|\theta|}_{c+\varepsilon
}-T^{|\theta|}_{c}
\bigr)
\mathop{\stackrel{(\mathit{law})}{\longrightarrow}}_{\varepsilon\rightarrow0}
\exp (2\beta_{T^{|\gamma
|}_{c}} ) T^{\gamma'}_{1},
\end{eqnarray}
where $\gamma'$ stands for a real Brownian motion, independent from
$\gamma$, and
\mbox{$T^{\gamma'}_{1}=\inf\{ t\dvt\gamma'_{t}=1 \}$}.
\end{longlist}
\end{prop}

\begin{pf}
We rely upon (\ref{skew-productplanar}) for the three proofs. By using
the scaling property of BM, we obtain:
\begin{eqnarray*}
T^{|\theta|}_{c}=A_{T^{|\gamma|}_{c}}(\beta) \stackrel{(\mathit{law})} {=}
A_{u}(\beta)|_{u=c^{2}T^{|\gamma|}_{1}}
\end{eqnarray*}
thus:
%
\begin{eqnarray}
\label{Tcbeta} \frac{1}{c^{2}} T^{|\theta|}_{c}\stackrel{(\mathit{law})}
{=}\int^{T^{|\gamma
|}_{1}}_{0} \mathrm{d}v   \exp (2c
\beta_{v} )  .
\end{eqnarray}
\begin{longlist}[(b)]
\item[(a)] For $c\rightarrow0$, the RHS of (\ref
{Tcbeta}) converges to $T^{|\gamma|}_{1}$,
thus we obtain part (a) of the proposition.\eject

\item[(b)] For $c\rightarrow\infty$, taking
logarithms on both sides of (\ref{Tcbeta}) and dividing by $c$, on the
LHS we obtain
$\frac{1}{c}  \log (T^{|\theta|}_{c} )-\frac{2}{c}\log
c$ and on the RHS:
\[
\frac{1}{c}   \log \biggl(\int^{T^{|\gamma|}_{1}}_{0}
\mathrm{d}v   \exp (2c\beta_{v} ) \biggr) = \log \biggl(\int
^{T^{|\gamma
|}_{1}}_{0} \mathrm{d}v   \exp (2c\beta_{v} )
\biggr)^{1/c},
\]
which, from the classical Laplace argument: $\|f\|_{p}\stackrel
{p\rightarrow\infty}{\longrightarrow}\|f\|_{\infty}$, converges for
$c\rightarrow\infty$, towards:
\[
2 \sup_{v\leq T^{|\gamma|}_{1}} (\beta_{v} ) \stackrel {(\mathit{law})} {=} 2|
\beta|_{T^{|\gamma|}_{1}}.
\]
\end{longlist}
This proves part (b) of the proposition.
\begin{longlist}[(c)]
\item[(c)]
%
\begin{eqnarray}
\label{Tthetatilde} T^{|\theta|}_{c+\varepsilon}-T^{|\theta|}_{c}&=&
\int^{T^{|\gamma
|}_{c+\varepsilon}}_{T^{|\gamma|}_{c}} \mathrm{d}u   \exp (2\beta_{u} )
\nonumber
\\
&=& \int^{T^{|\gamma|}_{c+\varepsilon}-T^{|\gamma|}_{c}}_{0} \mathrm{d}v   \exp (2
\beta_{T^{|\gamma|}_{c}} )   \exp \bigl(2 (\beta_{v+T^{|\gamma|}_{c}}-\beta_{T^{|\gamma|}_{c}}
) \bigr)
\\
&=& \exp (2\beta_{T^{|\gamma|}_{c}} ) \int^{T^{|\gamma
|}_{c+\varepsilon}-T^{|\gamma|}_{c}}_{0}
\mathrm{d}v   \exp (2B_{v} ),\nonumber
\end{eqnarray}
where $ (B_{s}\equiv\beta_{s+T^{|\gamma|}_{c}}-\beta_{T^{|\gamma|}_{c}},s\geq0 )$ is a BM independent of
$T^{|\gamma|}_{c}$.
\end{longlist}

We study now $\tilde{T}^{|\gamma|}_{c,c+\varepsilon}\equiv
T^{|\gamma|}_{c+\varepsilon}-T^{|\gamma|}_{c}$, the first hitting
time of the level $c+\varepsilon$ from $|\gamma|$, starting from $c$.
Thus, we define: $\rho_{u}\equiv|\gamma_{u}|$, starting also from
$c$. Thus, $\rho_{u}=c+\delta_{u}+L_{u}$, where $ (\delta_{s},s\geq0 )$ is a BM and $ (L_{s},s\geq0 )$ is the
local time of $\rho$ at 0.
Thus,
\begin{eqnarray}
\tilde{T}^{|\gamma|}_{c,c+\varepsilon}&=&\inf \{u\geq0\dvt \rho_{u}=c+
\varepsilon \}\equiv\inf \{u\geq0\dvt \delta_{u}+L_{u}=
\varepsilon \}
\nonumber\\[-8pt]\\[-8pt]
&\stackrel{u=\varepsilon^{2}v} {=}& \varepsilon^{2} \inf
\biggl\{ v\geq0\dvt \frac{1}{\varepsilon}\delta_{v\varepsilon^{2}}+\frac
{1}{\varepsilon}L_{v\varepsilon^{2}}=1
\biggr\} .\nonumber
\end{eqnarray}
From Skorokhod's lemma Revuz and Yor \cite{ReY99}:
\[
L_{u}=\mathop{\sup}_{y\leq u} \bigl( (-c-\delta_{y} )\vee0
\bigr)
\]
we deduce:
%
\begin{eqnarray}
\frac{1}{\varepsilon}L_{v\varepsilon^{2}}=\mathop{\sup}_{y\leq v\varepsilon^{2}}
\bigl( (-c-\delta_{y} )\vee 0 \bigr) \stackrel{y=
\varepsilon^{2}\sigma} {=} \mathop{\sup}_{\sigma\leq v} \biggl(
\biggl(-c-\varepsilon\frac{1}{\varepsilon}\delta_{\sigma
\varepsilon^{2}} \biggr)\vee0
\biggr)=0.
\end{eqnarray}
Hence, with $\gamma'$ denoting a new BM independent from $\gamma$,
(\ref{Tthetatilde}) writes:
%
\begin{eqnarray}
\label{Tthetatildebis} T^{|\theta|}_{c+\varepsilon}-T^{|\theta|}_{c}
&=& \exp (2\beta_{T^{|\gamma|}_{c}} ) \int^{\varepsilon
^{2}T^{\gamma'}_{1}}_{0} \mathrm{d}v
  \exp (2B_{v} ).
\end{eqnarray}
Thus, dividing both sides of $(\ref{Tthetatildebis})$ by $\varepsilon^{2}$ and making $\varepsilon\rightarrow0$, we obtain part (c) of the proposition.
\end{pf}

\begin{rem}\label{remc}
The asymptotic result (c) in Proposition \ref{prop1}
may also be obtained in a
straightforward manner from (\ref{Tthetatilde}) by
analytic computations. Indeed, using the Laplace transform of the first
hitting time of a
fixed level by the absolute value of a linear Brownian motion
$E [\mathrm{e}^{-{(\lambda^{2}}/{2})T^{|\gamma|}_{b}} ]=\frac
{1}{\cosh(\lambda b)}$
(see, e.g., Proposition 3.7, page 71 in Revuz and Yor \cite{ReY99}),
we have that for $0<c<b$, and $\lambda\geq0$:
%
\begin{eqnarray}
E \bigl[\mathrm{e}^{-({\lambda^{2}}/{2}) (T^{|\gamma|}_{b}-T^{|\gamma
|}_{c} )} \bigr] =\frac{\cosh(\lambda c)}{\cosh(\lambda b)}.
\end{eqnarray}
Using now $b=c+\varepsilon$, for every $\varepsilon>0$, the latter equals:
\[
\frac{\cosh({\lambda c}/{\varepsilon})}{\cosh (({\lambda}/{\varepsilon})(c+\varepsilon) )} \stackrel{\varepsilon\rightarrow0} {\longrightarrow}\mathrm{e}^{-\lambda}.
\]
The result follows now by remarking that $\mathrm{e}^{-\lambda}$ is the Laplace
transform (for the argument $\lambda^{2}/2$)
of the first hitting time of 1 by a linear Brownian motion $\gamma'$,
independent from $\gamma$.
\end{rem}
%

\subsection{Generalizations}

Obviously, we can obtain several variants of Proposition \ref{prop1},
by studying $T^{\theta}_{-bc,ac}$, $0<a,b\leq\infty$, for
$c\rightarrow0$ or $c\rightarrow\infty$,
and $a,b$ fixed. We define $T^{\gamma}_{-d,c}\equiv\inf\{ t\dvt\gamma_{t}\notin(-d,c) \}$ and we have:
\begin{itemize}
\item$\frac{1}{c^{2}}   T^{\theta}_{-bc,ac}
\mathop{\stackrel{(\mathit{law})}{\longrightarrow}}\limits_{c\rightarrow0}
 T^{\gamma}_{-b,a}$.
\item$\frac{1}{c}   \log (T^{\theta}_{-bc,ac} )
\mathop{\stackrel{(\mathit{law})}{\longrightarrow}}\limits_{c\rightarrow\infty}
2|\beta|_{T^{\gamma}_{-b,a}}$.
\end{itemize}
In particular, we can take $b=\infty$, hence the following corollary.

\begin{cor}\label{coro}
\begin{longlist}[(b)]
\item[(a)] For $c\rightarrow0$, we have
%
\begin{eqnarray}
\frac{1}{c^{2}}   T^{\theta}_{ac}
\mathop{\stackrel{(\mathit{law})}{\longrightarrow}}_{c\rightarrow0}
 T^{\gamma}_{a}.
\end{eqnarray}
\item[(b)] For $c\rightarrow\infty$, we have
%
\begin{eqnarray}
\label{asympclarge} \frac{1}{c}   \log \bigl(T^{\theta}_{ac}
\bigr)
\mathop{\stackrel{(\mathit{law})}{\longrightarrow}}_{c\rightarrow\infty}
 2|\beta
|_{T^{\gamma}_{a}} \stackrel{(\mathit{law})} {=}2|C_{a}|,
\end{eqnarray}
where $ (C_{a},a\geq0 )$ is a standard Cauchy process.
\end{longlist}
\end{cor}

\begin{rem}[(Yet another proof of Spitzer's theorem)]\label{rem}
Taking $a=1$, from Corollary \ref{coro}\textup{(b)}, we can obtain yet another
proof of Spitzer's celebrated asymptotic theorem stated in (\ref{Spi}).
Indeed, (\ref{asympclarge}) can be equivalently stated as:
%
\begin{eqnarray}
\label{asympP} P \bigl(\log T^{\theta}_{c}<cx \bigr)
\mathop{\stackrel{(\mathit{law})}{\longrightarrow}}_{c\rightarrow\infty}
P
\bigl(2|C_{1}|<x \bigr).
\end{eqnarray}
Now, the LHS of (\ref{asympP}) equals:
\begin{eqnarray}
P \bigl(\log T^{\theta}_{c}<cx \bigr)&\equiv& P
\bigl(T^{\theta
}_{c}<\exp(cx) \bigr)\equiv P \Bigl(
\sup_{u\leq\exp(cx)} \theta_{u} >c \Bigr)
\nonumber\\[-8pt]\\[-8pt]
&=& P \bigl(|\theta_{\exp(cx)}| >c \bigr)=P \biggl(|\theta_{t}| >
\frac{\log t}{x} \biggr),\nonumber
\end{eqnarray}
with $t=\exp(cx)$. Thus, because $|C_{1}|\stackrel
{{(\mathit{law})}}{=}|C_{1}|^{-1}$, (\ref{asympP}) now writes:
%
\begin{equation}
\mbox{for  every } x>0\mbox{ given}\qquad   P \biggl(|\theta_{t}| >\frac{\log
t}{x}
\biggr)
\mathop{\stackrel{(\mathit{law})}{\longrightarrow}}_{t\rightarrow\infty}
 P
\biggl(|C_{1}|>\frac{2}{x} \biggr),
\end{equation}
which yields precisely Spitzer's theorem (\ref{Spi}).
\end{rem}
%

\subsection{Speed of convergence}

We can easily improve upon Proposition \ref{prop1} by studying the
speed of convergence of the distribution of
$\frac{1}{c^{2}}   T^{|\theta|}_{c}$ towards that of $T^{|\gamma
|}_{1}$, that is, the following proposition.

\begin{prop}\label{propcvgspeed}
For any function $\varphi\in\mathcal{C}^{2}$, with compact support,
\begin{eqnarray}
\label{cvgspeed}
&& \frac{1}{c^{2}} \biggl(E \biggl[\varphi \biggl(
\frac
{1}{c^{2}}T^{|\theta|}_{c} \biggr) \biggr]-E \bigl[\varphi
\bigl(T^{|\gamma|}_{1} \bigr) \bigr] \biggr)
\nonumber\\[-8pt]\\[-8pt]
&&\quad                          \mathop{\longrightarrow}_{c\rightarrow0} E \biggl[\varphi' \bigl(T^{|\gamma
|}_{1}
\bigr) \bigl(T^{|\gamma|}_{1} \bigr)^{2} +
\frac{2}{3} \varphi'' \bigl(T^{|\gamma|}_{1}
\bigr) \bigl(T^{|\gamma|}_{1} \bigr)^{3} \biggr].\nonumber
\end{eqnarray}
\end{prop}
\begin{pf}
We develop $\exp (2c\beta_{v} )$, for $c\rightarrow0$, up
to the second order term, that is,
\begin{eqnarray*}
\mathrm{e}^{2c\beta_{v}}=1+2c\beta_{v}+2c^{2}\beta_{v}^{2}+
\cdots.
\end{eqnarray*}
More precisely, we develop up to the second order term, and we obtain
\begin{eqnarray*}
E \biggl[\varphi \biggl( \frac{1}{c^{2}}T^{|\theta|}_{c}
\biggr) \biggr] &=& E \biggl[\varphi \biggl( \int^{T^{|\gamma|}_{1}}_{0}
\mathrm{d}v   \exp (2c\beta_{v} ) \biggr) \biggr]
\\
&=& E \biggl[\varphi \bigl( T^{|\gamma|}_{1} \bigr)+
\varphi' \bigl( T^{|\gamma|}_{1} \bigr)  \int
^{T^{|\gamma|}_{1}}_{0} \bigl(2c\beta_{v}+2c^{2}
\beta_{v}^{2} \bigr)\,\mathrm{d}v \biggr]
\\
&&   {}  +\frac{1}{2} E \biggl[\varphi''
\bigl( T^{|\gamma
|}_{1} \bigr)  4c^{2} \biggl(\int
^{T^{|\gamma|}_{1}}_{0}\beta_{v}\,\mathrm{d}v
\biggr)^{2} \biggr] + c^{2}\mathrm{o}(c).
\end{eqnarray*}
We then remark that $E [\int^{t}_{0}\beta_{v}\,\mathrm{d}v ]=0$,
$E [\int^{t}_{0}\beta_{v}^{2}\,\mathrm{d}v ]=t^{2}/2$ and
$E [ (\int^{t}_{0}\beta_{v}\,\mathrm{d}v )^{2} ]=t^{3}/3$, thus
we obtain (\ref{cvgspeed}).
\end{pf}
%

\section{Checks via Bougerol's identity}\label{secchecks}

So far, we have not made use of Bougerol's identity (\ref{boug}), which
helps us to characterize the distribution of $T^{|\theta|}_{c}$ \cite
{Vak11}. In this subsection,
we verify that writing the Gauss--Laplace transform in (\ref{pseudoTL}) as:
%
\begin{eqnarray}
\label{pseudoTL2} E \biggl[ \sqrt{\frac{2}{\uppi}} \frac{1}{\sqrt{({1}/{c^{2}})T^{|\theta|}_{c}}} \exp
\biggl( -\frac{xc^{2}}{2T^{|\theta
|}_{c}} \biggr) \biggr] = \frac{1}{\sqrt{1+xc^{2}}}
\varphi_{m}\bigl(xc^{2}\bigr),
\end{eqnarray}
with $m=\uppi/(2c)$, we find asymptotically for $c\rightarrow0$ the
Gauss--Laplace transform of $T^{|\gamma|}_{1}$.
Indeed, from (\ref{pseudoTL2}), for $c\rightarrow0$, we obtain:
\begin{eqnarray}
\label{pseudoTL3} && E \biggl[ \sqrt{\frac{2}{\uppi}} \frac{1}{\sqrt{T^{|\gamma
|}_{1}}} \exp
\biggl( -\frac{x}{2T^{|\gamma|}_{1}} \biggr) \biggr]
\nonumber\\[-8pt]\\[-8pt]
&&   \quad            = \lim_{c\rightarrow0} \frac{2}{
(\sqrt{1+xc^{2}}+\sqrt{xc^{2}} )^{\uppi/(2c)}+ (\sqrt {1+xc^{2}}-\sqrt{xc^{2}} )^{\uppi/(2c)}}  .\nonumber
\end{eqnarray}
Let us now study:
\begin{eqnarray*}
\bigl(\sqrt{1+xc^{2}}+\sqrt{xc^{2}} \bigr)^{\uppi/(2c)}&=&
\exp \biggl(\frac{\uppi}{(2c)}\log \bigl[1+ \bigl(\sqrt{1+xc^{2}}-1
\bigr)+\sqrt {xc^{2}} \bigr] \biggr)
\\
&\thicksim& \exp \biggl(\frac{\uppi}{2c} \biggl[c\sqrt{x}+\frac
{xc^{2}}{2}
\biggr] \biggr)\mathop{\longrightarrow}_{c\rightarrow0}\exp \biggl(
\frac{\uppi\sqrt{x}}{2} \biggr).
\end{eqnarray*}
A similar calculation finally gives
%
\begin{eqnarray}
\label{pseudoTL4} E \biggl[ \sqrt{\frac{2}{\uppi}} \frac{1}{\sqrt{T^{|\gamma|}_{1}}} \exp
\biggl( -\frac{x}{2T^{|\gamma|}_{1}} \biggr) \biggr] = \frac{1}{\cosh (({\uppi}/{2})\sqrt{x} )}  ,
\end{eqnarray}
a result which is in agreement with the law of $\beta_{T^{|\gamma
|}_{1}}$, whose density is
%
\begin{eqnarray}
\label{betaTdensity} E \biggl[ \frac{1}{\sqrt{2\uppi T^{|\gamma|}_{1}}} \exp \biggl( -
\frac
{y^{2}}{2T^{|\gamma|}_{1}} \biggr) \biggr] = \frac{1}{2\cosh (({\uppi}/{2})y )}  .
\end{eqnarray}
Indeed, the law of $\beta_{T^{|\gamma|}_{c}}$ may be obtained from its
characteristic function which is given by Revuz and Yor \cite{ReY99}, page 73:
\[
E \bigl[ \exp(\mathrm{i}\lambda\beta_{T^{|\gamma|}_{c}}) \bigr] = \frac
{1}{\cosh(\lambda c)}  .
\]
It is well known that L\'{e}vy \cite{Lev80}, Biane and Yor \cite{BiY87}:
\begin{eqnarray}
\label{Fourier}
E \bigl[ \exp( \mathrm{i}\lambda\beta_{T^{|\gamma|}_{c}}) \bigr] &=&
\frac{1}{\cosh(\lambda c)} = \frac{1}{\cosh(\uppi\lambda{c}/{\uppi})} = \int^{\infty}_{-\infty}
\mathrm{e}^{\mathrm{i}  ( {\lambda c}/{\uppi}
) y } \frac{1}{2\uppi} \frac{1}{\cosh({y}/{2})}  \,\mathrm{d}y
\nonumber\\[-8pt]\\[-8pt]
&\stackrel{x={cy}/{\uppi}} {=}& \int^{\infty}_{-\infty}
\mathrm{e}^{\mathrm{i}\lambda x} \frac{1}{2\uppi} \frac{{\uppi}/{c}}{\cosh({x
\uppi}/{(2c)})} \, \mathrm{d}x = \int
^{\infty}_{-\infty} \mathrm{e}^{\mathrm{i} \lambda x} \frac{1}{2c}
\frac{1}{\cosh({x \uppi}/{(2c)})} \, \mathrm{d}x  .\qquad \nonumber
\end{eqnarray}
So, the density $h_{-c,c}$ of $\beta_{T^{|\gamma|}_{c}}$ is:\vspace*{-1pt}
\[
h_{-c,c}(y)= \biggl( \frac{1}{2c} \biggr) \frac{1}{\cosh({y\uppi}/{(2c)})} =
\biggl( \frac{1}{c} \biggr) \frac{1}{\mathrm{e}^{{y\uppi
}/{(2c)}} + \mathrm{e}^{-{y\uppi}/{(2c)}}}
\]
and for $c=1$, we obtain (\ref{betaTdensity}).

We recall from Remark \ref{remc} that (see also Pitman and Yor \cite{PiY03}, where
further results concerning the infinitely divisible distributions
generated by some L\'{e}vy processes associated with the hyperbolic
functions $\cosh$, $\sinh$ and $\tanh$ can also be found):
%
\begin{eqnarray}
E \biggl[ \exp \biggl( -\frac{\lambda^{2}}{2} T^{|\gamma|}_{c}
\biggr) \biggr] = \frac{1}{\cosh(\lambda c)}  ,
\end{eqnarray}
thus, for $c=1$ and $\lambda=\frac{\uppi}{2}\sqrt{x}$, (\ref
{pseudoTL4}) now writes:
%
\begin{eqnarray}
\label{pseudoTL5} E \biggl[ \sqrt{\frac{2}{\uppi}} \frac{1}{\sqrt{T^{|\gamma|}_{1}}} \exp
\biggl( -\frac{x}{2T^{|\gamma|}_{1}} \biggr) \biggr] = E \biggl[ \exp \biggl( -
\frac{x\uppi^{2}}{8} T^{|\gamma|}_{1} \biggr) \biggr],
\end{eqnarray}
a result which gives a probabilistic proof of the reciprocal relation
that was obtained in Biane, Pitman and Yor \cite{BPY01}
(using the notation of this article, Table 1, page 442):
\[
f_{C_{1}}(x)= \biggl(\frac{2}{\uppi x} \biggr)^{3/2}
f_{C_{1}} \biggl(\frac{4}{\uppi^{2} x} \biggr).\vspace*{-1pt}
\]

%

\printhistory

\end{document}